\magnification=1100
\input amstex
\documentstyle{amsppt}
\def\leftitem#1{\item{\hbox to\parindent{\enspace#1\hfill}}}

\def\dim{\operatorname{dim}}

\def\proclaim#1{\par\vskip.25cm\noindent{\bf#1. }\begingroup\it}
\def\endproclaim{\endgroup\par\vskip.25cm}
\leftheadtext{Shen, Wang} \rightheadtext{derivations of Lie algebras}

\topmatter
\title Derivations of Twisted Loop Algebras of Minimal Q-graded Subalgebras
\endtitle
\author  Yaxin Shen$^{\dagger}$, Xiandong Wang$^{\ddag}$
\endauthor
\affil   College of Mathematics and Statistics\\
       Qingdao University, Qingdao 266071\\
       P. R. China\\
       $^{\dag}$Email:  shenyaxin\@qdu.edu.cn \,\,\,\,\,\,\,\,\,\,\\
       $^{\ddag}$Email: wangxiandong\@qdu.edu.cn\\
\endaffil
\thanks Research supported by NSF of China (No.11472144).
\endthanks
 \abstract Derivations of twisted loop algebras of minimal Q-graded subalgebras of semisimple Lie algebras are investigated,
 and a decomposition formula of the whole derivation space is obtained.
 Homogenous almost inner derivations of twisted loop algebras and twisted affinizations of
 minimal Q-graded subalgebras are determined.
 \endabstract
\endtopmatter
\document

\smallskip\bigskip
\subhead 1. \ \ Introduction
\endsubhead
\medskip

Discussions on almost inner derivations of Lie algebras appear frequently in literature in the last several years [1-4,15].
One finds some classes of \,(mainly nilpotent and finite dimensional) Lie algebras in which almost inner derivations are determined.
In [5], the authors introduced a class of solvable Lie algebras called minimal Q-graded subalgebras of semisimple Lie algebras,
and checked that these Lie algebras have no non-inner almost inner derivations.
Moreover, almost inner derivations of (non-twisted) loop algebras and affinizations of these Lie algebras are also determined.

The concept of minimal Q-graded subalgebras of semisimple Lie algebras is motivated by Q-graded Lie algebras [6].
On the other hand, almost inner derivations of Lie algebras are related to almost inner automorphisms of Lie groups,
the latter is used to resolve isospectral problems of compact Riemannian manifolds. In [7-8], the authors
discussed on basic results concerning these problems.

In this paper, we provide new examples of Lie algebras in which derivations and homogenous almost inner derivations are determined.
We take a minimal Q-graded subalgebra $L$ as our starting point, and form the twisted loop algebra $\bar L$ and twisted affinization $\hat L$ for which
we investigate derivations and almost inner derivations. In section 2, we calculate derivations of a minimal Q-graded subalgebra $L$,
as bases for further investigation. In section 3, we deal with derivations of twisted loop algebras, and obtain
the following decomposition formula:
$$ Der(\bar L)=Der(L)_{\bar 0}\otimes R\oplus Der(L)_{\bar 1}\otimes Rt\oplus Cent(L)\otimes Der(S),$$
here $S=\Bbb F[t,t^{-1}]$ is the Laurent polynomial algebra and $R=\Bbb F[t^2,t^{-2}]$ is the subalgebra of $S$ consisting of even elements.

In section 4, we determine homogenous almost inner derivations of the twisted affinization $\hat L$ of a minimal Q-graded subalgebra $L$
with $\dim L=2\dim[L,L]$ just as we did in the untwisted case [5]:

Any homogenous element in the space $AID(\hat{L})/Inn(\hat{L})$ has the following form
$$ \sum_{1\leq i\leq l, j\in \Bbb{Z}}a_{ij}D_{ij}+Inn(\hat{L}),\,\, a_{ij}\in \Bbb{F}, $$
where the derivation $D_{ij}$ is an even almost inner derivation and determined by the following rules:
$$ D_{ij}([L,L]\otimes Rt)=0, \,\, D_{ij}(K)=0, $$
$$ D_{ij}(h_m\otimes t^{2n})=\delta_{im}\delta_{jn}K,\,\, 1\leq m\leq l, \,\, n\in \Bbb{Z},$$
here $h_1,\cdots ,h_l$ is a basis of the maximal torus $H$ which is contained in the subalgebra $L$.
Moreover, the infinite subset $\{D_{ij}+Inn(\hat{L}); 1\leq i\leq l, j\in \Bbb{Z}\}$ of the quotient space
$ AID(\hat{L})/Inn(\hat{L}) $ is linearly independent.

\vskip 5pt

\smallskip\bigskip
\subhead 2. \ \ Derivations of Minimal Q-graded Subalgebras
\endsubhead
\medskip

First we recall some notations from [5,9-10]: Let $\Bbb{F}$ be an algebraically closed field with characteristic $0$,
$\Bbb{L}$ a finite dimensional semisimple Lie algebra over $\Bbb{F}$,
$H$ a maximal torus, $\Phi$ the corresponding root system, and $\triangle=\{\alpha_{1},\alpha_{2},\cdots,\alpha_{l}\}$ a basis of $\Phi$,
hence we have the root space decomposition of $\Bbb{L}$:
$$ \Bbb{L}=H\oplus \bigoplus_{\alpha\in \Phi}L_{\alpha},$$
where $L_{\alpha}$ is the root space of $\Bbb{L}$ with respect to the root $\alpha$.

Let $Q=\Bbb{Z}{\Phi}$ be the free abelian group generated by $\triangle$,
a subalgebra $L$ of $\Bbb{L}$ is called a minimal Q-graded subalgebra if it satisfies the following conditions:

1) $H$ is contained in $L$, and there is a subset $\Psi$ of $\Phi$ such that $Q$ is spanned integrally by $\Psi$
in the sense that $Q=span_{\Bbb{Z}}{\Psi}$ and
$$ L=H\oplus\bigoplus_{\alpha\in{\Psi}}L_{\alpha}, \,\,\, \Psi{\subseteq}{\Phi}.$$

2) There is no subalgebra $L_{1}$ of $\Bbb{L}$ such that ${L_{1}}\subsetneq{L}$
which satisfies the same conditions as the subalgebra $L$ does as in 1).

Now we also suppose that the equality $\dim L=2\dim [L,L]$ holds on the Lie algebra $L$ as we did in sections 4 and 5 of our previous paper [5],
and so we can choose a basis of it like this:
$$ h_1,\cdots,h_l, \,\, x_1,\cdots, x_l,$$
here $\{h_1,\cdots,h_l\}$ is a basis of $H$ which is dual to the basis $\Psi=\{\beta_1,\cdots,\beta_l\}$ of the dual space $H^*$,
and $x_i\in L_{\beta_i}, 1\leq i\leq l$, these elements constitute a basis of the derived algebra $I=[L,L]$.

Note that the minimal Q-graded subalgebra $L$ just defined above is naturally a $\Bbb Z_2$-graded Lie algebra with
$L_{\bar 0}=H, L_{\bar 1}=[L,L]$, and we have the following decomposition of subspaces
$$ L=L_{\bar 0}\oplus L_{\bar 1}.$$
In order to determine all the derivations of the Lie algebra $L$, we first give a basic result on derivations of
$\Bbb Z_2$-graded Lie algebras.

\proclaim{Lemma 2.1} Let $\Cal L$ be any $\Bbb Z_2$-graded Lie algebra with the decomposition
$$ \Cal L={\Cal L}_{\bar 0}\oplus {\Cal L}_{\bar 1},$$
then we have the following decomposition of derivations of $\Cal L$:
$$ Der(\Cal L)=Der(\Cal L)_{\bar 0}\oplus Der(\Cal L)_{\bar 1},$$
where
$$ Der(\Cal L)_{\bar 0}=\{D\in Der(\Cal L); D(\Cal L_i)\subset {\Cal L}_i, i=\bar 0, \bar 1\},$$
$$ Der(\Cal L)_{\bar 1}=\{D\in Der(\Cal L); D(\Cal L_i)\subset {\Cal L}_{i+\bar 1}, i=\bar 0, \bar 1\}.$$

A derivation $D$ is called even if $D\in Der(\Cal L)_{\bar 0}$ or odd if $D\in Der(\Cal L)_{\bar 1}$,
even or odd derivations are usually called homogenous derivations.
\endproclaim
\demo{Proof} The decomposition formula in this lemma can be checked directly or by modifying the proof of Proposition 1.1 in [11].
\hfill $ $ \qed
\enddemo
\proclaim{Lemma 2.2} The minimal Q-graded subalgebra $L=H\oplus [L,L]$ is a $\Bbb{Z}_2$-graded Lie algebra with
$L_{\bar 0}=H, L_{\bar 1}=[L,L]$, and we have
$$ Der(L)=Der(L)_{\bar 0}\oplus Der(L)_{\bar 1} ,$$
where the two subspaces are characterized as follows:
$$ Der(L)_{\bar 0}=\{D\in Der(L); D(h_i)=0, D(x_i)=b_ix_i, b_i\in \Bbb{F}, 1\leq i\leq l\},$$
$$ Der(L)_{\bar 1}=\{D\in Der(L); D(h_i)=a_ix_i, D(x_i)=0, a_i\in \Bbb{F},1\leq i\leq l\}.$$
\demo{Proof} 1) Even case: Let $D\in Der(L)_{\bar 0}$, with basis $\{h_1,\cdots, h_l, x_1,\cdots, x_l\}$ of the Lie algebra $L$,
we have the following identities
$$ D(h_i)=\sum_k a_{ki}h_k, \,\, D(x_i)=\sum_k b_{ki}x_k, \,\, a_{ki}, b_{ki}\in \Bbb F,$$
$$ D[h_i,x_j]=[\sum_k a_{ki}h_k, x_j]+[h_i, \sum_k b_{kj}x_k]=a_{ji}x_j+b_{ij}x_i,$$
from which we deduce that $a_{ij}=b_{ij}=0$ whenever $i\neq j$ and $b_{ii}=a_{ii}+b_{ii}$ for any index $i$,
hence the derivation $D$ acts on the basis as required.

2) Odd case:  Let $D\in Der(L)_{\bar 1}$, with the same basis of the Lie algebra $L$ as in the even case,
we have the following identities
$$ D(h_i)=\sum_k a_{ki}x_k, \,\, D(x_i)=\sum_k c_{ki}h_k, \,\, a_{ki}, c_{ki}\in \Bbb F,$$
$$ D[h_i,h_j]=[\sum_k a_{ki}x_k, h_j]+[h_i, \sum_k a_{kj}x_k]=-a_{ji}x_j+a_{ij}x_i,$$
from which we deduce that $a_{ij}=0$ whenever $i\neq j$, and $c_{ki}=0 \,(\forall k,i) $ as $x_i\in [L,L]$ for any $i$,
the derivation $D$ acts on the given basis as required.

Conversely, in both of the two cases above, it can be checked easily that a linear transformation $D\in End L$ defined accordingly is an even or odd
derivation of the Lie algebra $L$, the proof of the lemma is finished.
\hfill $ $ \qed
\enddemo

\proclaim{Corollary 2.3} With notations as in Lemma 2.2, we have the following identity which means that
any derivation of the Lie algebra $L$ is an inner derivation:
$$  Der(L)=Inn(L).$$
\endproclaim
\demo{Proof} Suppose that $D$ is an inner derivation of the Lie algebra $L$, and there exists an element $y=\sum_kb_kh_k-\sum_ka_kx_k\in L$ such that
$ D=ad y$, then we have
$$ D(h_i)=[\sum_kb_kh_k-\sum_ka_kx_k, h_i]=a_ix_i,$$
$$ D(x_i)=[\sum_kb_kh_k-\sum_ka_kx_k, x_i]=b_ix_i.$$

Now it is obvious from descriptions in Lemma 2.2 and identities above that inner derivations of the Lie algebra $L$
exhaust all derivations of $L$, hence any derivation of $L$ is inner, and the proof of the corollary is finished.
\hfill  $ $ \qed
\enddemo

Remark: From the concrete actions of general derivations of $L$ in Lemma 2.2,
it is easy to see that any local derivation of the Lie algebra $L$ is a derivation.
Roughly speaking, a local derivation of a Lie algebra $L$ is a linear transformation of $L$ that looks like a derivation of $L$ locally.
For more information about local derivations of Lie algebras, see references [13-14].

\vskip 5pt

\smallskip\bigskip
\subhead 3. \ \ Derivations of Twisted Loop Algebras
\endsubhead
\medskip

In this section, we investigate derivations of the twisted loop algebra $\bar L$ of a minimal Q-graded subalgebra $L$,
using the natural ${\Bbb Z}_2$-gradation of the derivation space $Der(\bar L)$ of the Lie algebra $\bar L$.

\proclaim{Definition 3.1} With notations as before, define the twisted loop algebra $\bar L$ and twisted affinization $\hat{L}$ of
a minimal Q-graded subalgebra $L$ as follows
$$ \bar L=H\otimes \Bbb{F}[t^2,t^{-2}]\oplus I\otimes \Bbb{F}[t^2,t^{-2}]t,$$
$$ \hat{L}=H\otimes \Bbb{F}[t^2,t^{-2}]\oplus I\otimes \Bbb{F}[t^2,t^{-2}]t\oplus \Bbb{F}K. $$
Hence the brackets of the Lie algebras $\bar L$ and $\hat L$ are induced by the brackets of the loop algebra $L\otimes \Bbb{F}[t,t^{-1}]$
and affinization $\tilde{L}=L\otimes \Bbb{F}[t,t^{-1}]\oplus \Bbb FK$ of $L$ respectively
 (for brevity, we will denote $\Bbb{F}[t^2,t^{-2}]$ by $R$ and $\Bbb{F}[t,t^{-1}]$ by $S$ from now on).
 \endproclaim

Set $\bar{L}_{\bar 0}=H\otimes R$ and $\bar{L}_{\bar 1}=I\otimes Rt$,
then we see that $\bar L=\bar{L}_{\bar 0}\oplus \bar{L}_{\bar 1}$ is a $\Bbb{Z}_2$-graded Lie algebra,
and the subspaces $\bar{L}_{\bar 0}$ and $\bar{L}_{\bar 1}$ have respectively the following bases:
$$ h_i\otimes t^{2j},  \,\, 1\leq i\leq l, j\in \Bbb{Z}; $$
$$ x_i\otimes t^{2j+1}, \,\, 1\leq i\leq l, j\in \Bbb{Z}.$$
Similarly, set $\hat{L}_{\bar 0}=H\otimes R\oplus \Bbb{F}K$ and $\hat{L}_{\bar 1}=I\otimes Rt$,
then $\hat L=\hat{L}_{\bar 0}\oplus \hat{L}_{\bar 1}$ is a $\Bbb{Z}_2$-graded Lie algebra,
and the subspaces $\hat{L}_{\bar 0}$ and $\hat{L}_{\bar 1}$ have respectively the following bases:
$$ h_i\otimes t^{2j}, \, K,  \,\, 1\leq i\leq l, j\in \Bbb{Z}; $$
$$ x_i\otimes t^{2j+1}, \,\, 1\leq i\leq l, j\in \Bbb{Z}.$$

Note that the vector spaces $\bar{L}_{\bar 0}$, $\bar{L}_{\bar 1}$ and $\bar L$ are natural $R$-$R$-bimodules.

\proclaim{Lemma 3.2} With notations as above, we have the following decomposition of vector subspaces:
$$ Der(\bar{L})=Der(\bar L)_{\bar 0}\oplus Der(\bar L)_{\bar 1},$$
where
$$ Der(\bar L)_{\bar 0}=\{D\in Der(\bar L); D(\bar{L}_i)\subset \bar{L}_i, i=\bar 0, \bar 1\},$$
$$ Der(\bar L)_{\bar 1}=\{D\in Der(\bar L); D(\bar{L}_i)\subset \bar{L}_{i+\bar 1}, i=\bar 0, \bar 1\}.$$
\endproclaim
\demo{Proof} This is a special case of Lemma 2.1.
\hfill $ $ \qed
\enddemo

It is obvious that $\bar L$ is isomorphic to the quotient algebra $\hat{L}/\Bbb{F}K$,
so we first describe derivations of the Lie algebra $\bar L$. The following result on derivations is similar
to Lemma 4.1 in [5] or Lemma 2.3 in [12].

\proclaim{Lemma 3.3} With notations as above, we have a decomposition formula below
$$ Der(\bar{L})=D_R(\bar{L})\oplus D_{L\otimes 1}(\bar{L}),$$
where the subspace $D_R(\bar L)$ consists of derivations $d\in Der(\bar L)$ satisfying
$$ d(h\otimes s_1s_2)=s_1d(h\otimes s_2)=d(h\otimes s_1)s_2, \forall h\in H,\forall s_1,s_2\in R,$$
$$ d(x\otimes s_1s_2t)=s_1d(x\otimes s_2t)=d(x\otimes s_1t)s_2, \forall x\in I,\forall s_1,s_2\in R.$$
and the subspace $D_{L\otimes 1}(\bar L)$ consists of derivations $\delta\in Der(\bar L)$ satisfying
$$ \delta(H\otimes 1)=0, \,\, \delta(I\otimes t)=0.$$
\endproclaim
\demo{Proof} For any derivation $D\in Der(\bar{L})$, we define an element $d\in D_R(\bar{L})$, this is a linear map $d\in End\bar{L}$
determined by:
$$ d(h\otimes s)=D(h\otimes 1)s, \,\, d(x\otimes st)=D(x\otimes t)s, \,\, h\in H, x\in I, s\in R. $$

To check that $d$ is a derivation of the Lie algebra $\bar{L}$, we calculate as follows: for any elements $h,h'\in H$ and $s,s'\in R$, we have
$$ [d(h\otimes s), h'\otimes s']=[D(h\otimes 1)s,h'\otimes s']=[D(h\otimes 1), h'\otimes 1]ss',$$
$$ [h\otimes s, d(h'\otimes s')]=[h\otimes s,D(h'\otimes 1)s']=[h\otimes 1, D(h'\otimes 1)]ss'.$$
So, by adding the two sides of the identities above we obtain
$$ [d(h\otimes s), h'\otimes s']+[h\otimes s, d(h'\otimes s')]=0=d[h\otimes s, h'\otimes s'].$$
Similarly, we have
$$ [d(h\otimes s), x\otimes s't]+[h\otimes s, d(x\otimes s't)]=d[h\otimes s, x\otimes s't],$$
$$ [d(x\otimes st), y\otimes s't]+[x\otimes st, d(y\otimes s't)]=d[x\otimes st, y\otimes s't]$$
from which we see that $d$ is a derivation, and $d\in D_R(\bar{L})$.

Finally, let $\delta=D-d$, then it is clear that $\delta(h\otimes 1)=\delta(x\otimes t)=0$ for any elements $h\in H$ and $x\in I$.
Hence the required decomposition holds.
\hfill $ $ \qed
\enddemo

\proclaim{Lemma 3.4} Notations as before, we have the following identity:
$$ D_R(\bar L)=Der(L)_{\bar 0}\otimes R\oplus Der(L)_{\bar 1}\otimes Rt.$$

In particular, any homogenous almost inner derivation $D\in D_{R}(\bar L)$ is an inner derivation of the twisted loop algebra $\bar L$.
\endproclaim
\demo{Proof} 1) Even case: Let $D$ be an even derivation in $D_R(\bar L)$, it is determined by its action
on elements of the form $h\otimes 1, h\in H$, $x\otimes t, x\in I$. Since the Lie algebra $L$ is finite dimensional,
there exist integers $r_1, r_2$ such that the following identities hold:

$$ D(h\otimes 1)=\sum_{i=r_1}^{r_2}D_i(h)\otimes t^{2i},  \,\, h\in H, $$
$$ D(x\otimes t)=\sum_{i=r_1}^{r_2}D_i(x)\otimes t^{2i+1}, \,\, x\in I.$$

It can be checked that the linear map $D_i\in End L$ is an even derivation of the Lie algebra $L$, $r_1\leq i\leq r_2$.
For example, we have
$$ D([h,x]\otimes t)=D[h\otimes 1, x\otimes t]=[D(h\otimes 1), x\otimes t]+[h\otimes 1, D(x\otimes t)] $$
$$ =\sum [D_i(h)\otimes t^{2i}, x\otimes t]+\sum [h\otimes 1, D_i(x)\otimes t^{2i+1}] $$
$$ =\sum [D_i(h), x]\otimes t^{2i+1}+\sum [h,D_i(x)]\otimes t^{2i+1}$$
which means that
$$ D_i[h,x]=[D_i(h), x]+[h,D_i(x)], \,\, \forall h\in H,\,\, \forall x\in I.$$

To show that $D\in Der(L)_{\bar 0}\otimes R$, we calculate its action on general elements: For any elements $h\in H, x\in I, s\in R$, we have
$$ D(h\otimes s)=\sum D_i(h)\otimes t^{2i}s=(\sum D_i\otimes t^{2i})(h\otimes s),$$
$$ D(x\otimes st)=\sum D_i(x)\otimes t^{2i+1}s=(\sum D_i\otimes t^{2i})(x\otimes st)$$
from which we deduce that $D=\sum D_i\otimes t^{2i}\in Der(L)_{\bar 0}\otimes R$. Conversely, it is easy to see that any element
in $Der(L)_{\bar 0}\otimes R$ is an even derivation in $D_R(\bar L)$.

If additionally $D$ is almost inner, for any elements $h\in H, x\in I$, there exist some elements $a_k,c_k\in H, b_k,d_k\in I$ \,$(1\leq k\leq l)$
satisfying
$$ D(h\otimes 1)=[h\otimes 1, \sum a_k\otimes t^{2k}+\sum b_k\otimes t^{2k+1}]=0,$$
$$ D(x\otimes t)=[x\otimes t, \sum c_k\otimes t^{2k}+\sum d_k\otimes t^{2k+1}] = \sum [x,c_k]\otimes t^{2k+1}$$
from which we deduce that $D_i(h)=0, D_i(x)=[x, c_i]$ and $D_i$ is an almost inner derivation of the Lie algebra $L$.
Hence $D_i$ is inner by Theorem 3.5 in [5].

Let $D_i=ad y_i, r_1\leq i\leq r_2$, then $D_i(h)=[y_i, h]=0$ for any $h\in H$, so $y_i\in H$, and the following identities
imply that $D$ is an inner derivation:
$$ \eqalign{ &
   D(h\otimes 1)=\sum D_i(h)\otimes t^{2i}=0, \cr
  & D(x\otimes t)=\sum ad y_i(x)\otimes t^{2i+1}=[\sum y_i\otimes t^{2i}, x\otimes t], \cr
  & D(\sum a_k\otimes t^{2k}+\sum b_k\otimes t^{2k+1}])=\sum D(a_k\otimes 1)t^{2k}+\sum D(b_k\otimes t)t^{2k}\cr
  & = \sum_k\sum_i ad y_i(b_k)\otimes t^{2i+2k+1}=[\sum y_i\otimes t^{2i}, \sum b_k\otimes t^{2k+1}]\cr
  & = [\sum y_i\otimes t^{2i},\sum a_k\otimes t^{2k}+\sum b_k\otimes t^{2k+1}],\cr }$$
where we take arbitrary elements $h\in H, x\in I, a_k\in H, b_k\in I, 1\leq k\leq l$.

2) Odd case: Similarly,
let $D$ be an odd derivation in $D_R(\bar L)$, it is determined by its action
on elements of the form $h\otimes 1, h\in H$, $x\otimes t, x\in I$, and
there exist integers $r_1, r_2$ such that the following identities hold:

$$ D(h\otimes 1)=\sum_{i=r_1}^{r_2}D_i(h)\otimes t^{2i-1},  \,\, h\in H, $$
$$ D(x\otimes t)=\sum_{i=r_1}^{r_2}D_i(x)\otimes t^{2i}, \,\, x\in I.$$

It can be checked that the linear map $D_i\in End L$ is an odd derivation of the Lie algebra $L$, $r_1\leq i\leq r_2$.
For example, we have
$$ D([h,x]\otimes t)=D[h\otimes 1, x\otimes t]=[D(h\otimes 1), x\otimes t]+[h\otimes 1, D(x\otimes t)] $$
$$ =\sum [D_i(h)\otimes t^{2i-1}, x\otimes t]+\sum [h\otimes 1, D_i(x)\otimes t^{2i}] $$
$$ =\sum [D_i(h), x]\otimes t^{2i}+\sum [h,D_i(x)]\otimes t^{2i}$$
which means that
$$ D_i[h,x]=[D_i(h), x]+[h,D_i(x)], \,\, \forall h\in H,\,\, \forall x\in I.$$

To show that $D\in Der(L)_{\bar 1}\otimes Rt$, we calculate its action on general elements: For any elements $h\in H, x\in I, s\in R$, we have
$$ D(h\otimes s)=\sum D_i(h)\otimes t^{2i-1}s=(\sum D_i\otimes t^{2i-1})(h\otimes s),$$
$$ D(x\otimes st)=\sum D_i(x)\otimes t^{2i}s=(\sum D_i\otimes t^{2i-1})(x\otimes st)$$
from which we deduce that $D=\sum D_i\otimes t^{2i-1}\in Der(L)_{\bar 1}\otimes Rt$. Conversely, it is easy to see that any element
in $Der(L)_{\bar 1}\otimes Rt$ is an odd derivation in $D_R(\bar L)$.

If additionally $D$ is almost inner, we can calculate brackets as just we did before in the even case,
and conclude that $D$ is an inner derivation of $\bar L$.

3) General case: Suppose that $D\in D_R(\bar L)$ and write it as a sum: $D=D_{\bar 0}+D_{\bar 1}$,
 where $D_{\bar 0}\in Der(\bar L)_{\bar 0}, D_{\bar 1}\in Der(\bar L)_{\bar 1}$, we will show below that $D_{\bar 0}\in D_R(\bar L)$
 which together with 1)-2) above imply the following relation:
 $$ D_R(\bar L)\subset Der(L)_{\bar 0}\otimes R\oplus Der(L)_{\bar 1}\otimes Rt.$$

For arbitrary elements $h\in H, x\in I, s\in R$, the following identities insure that the derivation $D_{\bar 0}\in D_R(\bar L)$,
hence the above relation holds:
$$ D(h\otimes s)=D(h\otimes 1)s, \,\, \,\, D(x\otimes st)=D(x\otimes t)s,$$
$$ D_{\bar 0}(h\otimes s)+D_{\bar 1}(h\otimes s)=D_{\bar 0}(h\otimes 1)s+D_{\bar 1}(h\otimes 1)s,$$
$$ D_{\bar 0}(x\otimes st)+D_{\bar 1}(x\otimes st)=D_{\bar 0}(x\otimes t)s+D_{\bar 1}(x\otimes t)s,$$
$$ D_{\bar 0}(h\otimes s)=D_{\bar 0}(h\otimes 1)s, \,\,\,\, D_{\bar 0}(x\otimes st)=D_{\bar 0}(x\otimes t)s. $$

Since the opposite direction of the relation can be checked easily, we conclude that the decomposition formula of the lemma holds.
\hfill $ $ \qed
\enddemo

\proclaim{Lemma 3.5} Suppose that $D\in D_{L\otimes 1}(\bar L)$, and $h_1,\cdots,h_l, x_1,\cdots, x_l$
is the basis of $L$ defined as before, then $D$ is an even derivation and it is determined
by some Laurent polynomials $f_{11}, \cdots, f_{l1}\in R$ such that
$$ D(h_i\otimes t^{2j})= h_i\otimes f_{ij}(t), \,\, D(x_i\otimes t^{2j+1})=x_i\otimes f_{ij}(t)t,$$
$$ f_{ij}(t)=jf_{i1}(t)t^{2j-2},\,\, 1\leq i\leq l, \, j\in \Bbb{Z},$$

In particular, any homogenous almost inner derivation $D\in D_{L\otimes 1}(\bar L)$ of the twisted loop algebra $\bar L$ is zero.
\endproclaim
\demo{Proof} 1) Even case: Let $D$ be an even derivation in $D_{L\otimes 1}(\bar L)$,
by definition we have $D(h\otimes 1)=D(x\otimes t)=0$, for any $h\in H, x\in I$.
Also the derivation $D$ is determined by its action on basis elements:
$$ D(h_i\otimes t^{2j})=\sum h_k\otimes f_{ij}^k(t), $$
$$ D(x_i\otimes t^{2j+1})=\sum x_k\otimes g_{ij}^k(t)t,$$
where $f_{ij}^k(t)$ and $g_{ij}^k(t)$ are Laurent polynomials in $R$, $1\leq i,k\leq l, j\in \Bbb{Z}$.

From the following identities
$$ \eqalign{ & D[h_i\otimes t^{2j},x_m\otimes t^{2n+1}] \cr
             & =[\sum h_k\otimes f_{ij}^k(t), x_m\otimes t^{2n+1}]+ [h_i\otimes t^{2j},\sum x_k\otimes g_{mn}^k(t)t] \cr
             & = x_m\otimes t^{2n+1}f_{ij}^m(t)+x_i\otimes t^{2j+1}g_{mn}^i(t),\cr
             & = \delta_{im}D(x_m\otimes t^{2j+2n+1}) \cr
             & = \delta_{im}\sum x_k\otimes g_{m,j+n}^k(t)t, \cr}$$
we deduce that $f_{ij}^m(t)=g_{mn}^i(t)=0$ whenever $i\neq m$. For $i=m$, we have
$$ g_{m,j+n}^m(t)=f_{mj}^m(t)t^{2n}+g_{mn}^m(t)t^{2j},$$
this implies $f_{ij}^i(t)=g_{ij}^i(t)$ since the identity $g_{i0}^i(t)=0$ holds for any $i$.
We thus obtain the simplified action on basis elements
$$ D(h_i\otimes t^{2j})= h_i\otimes f_{ij}(t), \,\, D(x_i\otimes t^{2j+1})=x_i\otimes f_{ij}(t)t,$$
$$ f_{i,j+1}(t)=f_{ij}(t)t^2+f_{i1}(t)t^{2j},\,\, 1\leq i\leq l, \, j\in \Bbb{Z},$$
or
$$ f_{ij}(t)=jf_{i1}(t)t^{2j-2}, \,\, 1\leq i\leq l, \, j\in \Bbb{Z}, $$
where $f_{11}(t),\cdots , f_{l1}(t)$ are Laurent polynomials in $R$,  which determine uniquely the even derivation $D\in D_{L\otimes 1}(\bar L)$.

In particular, if $D$ is an almost inner derivation, then $D(\bar L)\subset [\bar L, \bar L]$, hence we get that $f_{ij}(t)=0, \forall i,j$,
and thus $D$ is zero.

2) Odd case: Let $D\in D_{L\otimes 1}(\bar L)$ be an odd derivation of the Lie algebra $\bar L$,
so that $D(h\otimes 1)=D(x\otimes t)=0$, for any $h\in H, x\in I$.
Also we may suppose that
$$ D(h_i\otimes t^{2j})=\sum x_k\otimes f_{ij}^k(t)t, $$
$$ D(x_i\otimes t^{2j+1})=\sum h_k\otimes g_{ij}^k(t),$$
where $f_{ij}^k(t)$ and $g_{ij}^k(t)$ are Laurent polynomials in $R$, $1\leq i,k\leq l, j\in \Bbb{Z}$.
We have the following identities:
$$ \eqalign{ & 0=D[h_i\otimes t^{2j}, h_m\otimes t^{2n}] \cr
             & = [\sum_k x_k\otimes f_{ij}^k(t)t, h_m\otimes t^{2n}]+[h_i\otimes t^{2j},\sum_k x_k\otimes f_{mn}^k(t)t]\cr
             & = -x_m\otimes f_{ij}^m(t) t^{2n+1}+ x_i\otimes f_{mn}^i(t) t^{2j+1} \cr}$$
from which we deduce that $f_{ij}^m(t)=0$ whenever $m\neq i$, and $f_{mj}^m(t)t^{2n}=f_{mn}^mt^{2j}$
for any integers $n, j\in \Bbb{Z}, 1\leq m\leq l$.

Since $f_{i0}^i(t)=0, \forall i$, we deduce that $f_{ij}^k(t)=0, \forall i,j,k$,
and $D(h_i\otimes t^{2j})=0 $ for any index $1\leq i\leq l, j\in \Bbb{Z}$.
Similarly, $D(x_i\otimes t^{2j+1})=0$ for any $1\leq i\leq l, j\in \Bbb{Z}$, therefore the derivation $D$ is zero.

3) General case: Let $D\in D_{L\otimes 1}(\bar L)$, and $D=D_{\bar 0}+D_{\bar 1}$, $D_{\bar 0}$ is even and $D_{\bar 1}$ is odd,
then obviously $D_{\bar 0}$ and $D_{\bar 1}$ belong to $D_{L\otimes 1}(\bar L)$, hence $D_{\bar 1}=0$ by 2) and the proof of the lemma is finished.
\hfill $ $ \qed
\enddemo

Recall from [5] that the centroid $Cent(L)$ of the minimal Q-graded subalgebra $L$ is commutative, and a basis
$ \lambda_1,\cdots, \lambda_l$ of $Cent(L)$ can be chosen satisfying
$$ \lambda_i(h_j)=\delta_{ij}h_j, \,\, \lambda_i(x_j)=\delta_{ij}x_j, \,\, 1\leq i,j\leq l.$$
Using $\lambda_1,\cdots, \lambda_l$, the derivation $D$ in Lemma 3.5 can be rewritten as follows:
$$ D(h_i\otimes t^{2j})= h_i\otimes jf_{i1}(t)t^{2j-2}=\sum \lambda_k(h_i)\otimes jt^{2j-2}f_{k1}(t),$$
$$ D(x_i\otimes t^{2j+1})= x_i\otimes jf_{i1}(t)t^{2j-1}=\sum \lambda_k(x_i)\otimes jt^{2j-2}f_{k1}(t)t.$$

The following theorem can be regarded as an analogous result to Theorem 4.8 in [5] in the untwisted case.

\proclaim{Theorem 3.6} Suppose that $\bar L$ is the twisted loop algebra of the minimal Q-graded subalgebra $L$ with
$\dim L=2\dim[L,L]$, then the derivation algebra $Der(\bar L)$ has a decomposition
$$ Der(\bar L)=Der(L)_{\bar 0}\otimes R\oplus Der(L)_{\bar 1}\otimes Rt\oplus Cent(L)\otimes Der(S),$$
where we identify an element $d\in Der(S)$ with a linear map $\delta: R\rightarrow R$ such that
$$ \delta(t^{2j})=jt^{2j-2}f(t),\,\, \forall j\in \Bbb{Z},$$
and $f(t)\in R$ is a Laurent polynomial depended only on $d$.

Moreover any homogenous almost inner derivation of the Lie algebra $\bar L$ is an inner derivation.
\endproclaim
\demo{Proof} From Lemma 3.3-3.5, the decomposition formula follows. We need only to prove the second statement of the theorem.

According to Lemma 3.3, any homogenous almost inner derivation $D$ can be written as follows
$$ D=D_1+D_2, \,\, D_1\in D_R(\bar L), D_2\in D_{L\otimes 1}(\bar L), $$
where $D_1$ and $D_2$ are homogenous derivations having the same parity with $D$ by the proof of Lemma 3.3.
Applying Lemma 3.4 and 3.5, it suffices to show that the homogenous derivation $D_1$ is an almost inner derivation.

1) Suppose that $D_1$ is even: Since $D$ is an almost inner derivation,
for any elements $h\in H, s\in R$, there exist elements $f_i, \varphi_i, \psi_{ik}\in R$ \,($1\leq i,k\leq l$) such that
$$ D_1(h\otimes s)=D(h\otimes 1)s=[h\otimes 1, \sum_k h_k\otimes f_k]s=0, $$
$$ D(x_i\otimes t)=[x_i\otimes t, \sum_k h_k\otimes \psi_{ik}]= -x_i\otimes \psi_{ii}t=[x_i\otimes t, \sum_k h_k\otimes \psi_{kk}],$$
$$ D_1(\sum_i x_i\otimes \varphi_i t)=\sum_i D(x_i\otimes t)\varphi_i$$
$$ =\sum_i[x_i\otimes t, \sum_k h_k\otimes \psi_{kk}]\varphi_i =[\sum_i x_i\otimes \varphi_it, \sum_k h_k\otimes \psi_{kk}],$$
here $h_1,\cdots,h_l, x_1,\cdots, x_l$ is the basis of the Lie algebra $L$ as defined in Lemma 3.5.
From these formulas, we confirm that the derivation $D_1$ is indeed almost inner.

2) Suppose that $D_1$ is odd: With the same notations as in 1) and $x\in I$,
using that $H$ and $I$ are abelian, we have the following identities
$$ D_1(x\otimes st)=D(x\otimes t)s=[x\otimes t, \sum_k x_k\otimes f_k]s=0, $$
$$ D(h_i\otimes 1)=[h_i\otimes 1, \sum_k x_k\otimes \psi_{ik}]= x_i\otimes \psi_{ii}=[h_i\otimes 1, \sum_k x_k\otimes \psi_{kk}],$$
$$ D_1(\sum_i h_i\otimes \varphi_i )=\sum_i D(h_i\otimes 1)\varphi_i$$
$$ =\sum_i[h_i\otimes 1, \sum_k x_k\otimes \psi_{kk}]\varphi_i =[\sum_i h_i\otimes \varphi_i, \sum_k x_k\otimes \psi_{kk}] $$
from which we deduce easily that $D_1$ is an almost inner derivation as required.
\hfill $ $ \qed
\enddemo

\vskip 5pt

\smallskip\bigskip
\subhead 4. \ \ Almost Inner Derivations of Twisted Affinizations
\endsubhead
\medskip

Now we are ready to describe almost inner derivations of the twisted affinization $\hat L$ of the minimal Q-graded subalgebra $L$ of
a semisimple Lie algebra.

\proclaim{Lemma 4.1} Suppose that $D\in Der(\hat L)$ is a homogenous almost inner derivation, then there exists some
$y\in \hat L$ such that the following relations hold
$$ (D-ad y)(I\otimes Rt)=0, \,\, (D-ad y)(K)=0, $$
$$ (D-ad y)(H\otimes R)\subset \Bbb{F}K. $$

\endproclaim
\demo{Proof} By assumption the derivation $D$ is almost inner, we have $D(K)\in [K, \hat{L}]$ and so $D(K)=0$.
Thus the derivation $D$ induces a map $\bar{D}$ as follows
$$ \bar{D}: \hat{L}/\Bbb{F}K\rightarrow \hat{L}/\Bbb{F}K, \,\, \bar{w}\mapsto \overline{D(w)}.$$
It can be checked that the map $\bar{D}$ is a homogenous almost inner derivation of the quotient algebra $\hat{L}/\Bbb{F}K$,
which is isomorphic to the twisted loop algebra $\bar L$.

From Theorem 3.6, we deduce that there exists some element $y\in \hat{L}$ such that $\bar{D}(\bar{w})=ad(\bar{y})(\bar{w}), \forall w\in \hat{L}$,
and therefore we have $(D-ad(y))(\hat{L})\subset \Bbb{F}K$.
In particular, we obtain that $(D-ad(y))[\hat{L}, \hat{L}]=0 $ and all the relations of the lemma hold for this element $y$.
\hfill $ $ \qed
\enddemo

\proclaim{Lemma 4.2} Notations as defined previously, for any integers $i,j\in \Bbb{Z}$ with $1\leq i\leq l$,
define a linear map $D_{ij}$ by
$$ D_{ij}: \hat L\rightarrow \hat L, \,\, K\mapsto 0,$$
$$ x_m\otimes t^{2n+1}\mapsto 0, \,\, h_m\otimes t^{2n}\mapsto \delta_{im}\delta_{jn}K, \,\, 1\leq m\leq l, n\in \Bbb{Z},$$
then $D_{ij}$ is an even almost inner derivation of the affinization $\hat L$.
\endproclaim
\demo{Proof} Take a basis $\{h_1,\cdots, h_l, x_1,\cdots, x_l\}$ of $L$ as before: $\{h_1,\cdots, h_l\}$ is a basis of $H$ which is dual to
the basis $\Psi=\{\beta_1,\cdots,\beta_l\}$ of $H^*$, $\{x_1,\cdots, x_l\}$ is a basis of $I=[L,L]$ with $x_i\in L_{\beta_i}, \forall i$,
so we have $[h_i, x_j]=\delta_{ij}x_j$. Also we have another basis $\{h_1',\cdots, h_l'\}$ of $H$ dual to $\{h_1,\cdots,h_l\}$
with respect to the bilinear form of $L$,
and we have $(h_i,h_j')=\delta_{ij}, \forall i,j$.

For any element $X=\sum b_{mn}h_m\otimes t^{2n}+\sum c_{uv}x_u\otimes t^{2v+1}\in \hat L$ without central component\,
(this is enough for the proof of the lemma),
we will find an element
$$ Y=\frac{1}{2}j^{-1}h_i'\otimes t^{-2j}-\sum_{k\in B} d_kh_k'\otimes t^{-2j}+\sum e_{pq}x_p\otimes t^{2q+1}$$
satisfying that $D_{ij}(X)=[X,Y]$, here we set
$$ A=\{m\in \{1,\cdots, l\}; \exists n, s.t.,  b_{mn}\neq 0\}, \,\, B=\{1,\cdots, l\}\backslash A.$$
Regarding these coefficients $d_k$ and $e_{pq}$ as unknowns,
we will solve some equations to get their values for the existence of the element $Y$,
this can be done by the calculations below
$$ \eqalign{ & b_{ij}K  =D_{ij}(X)=[X,Y]
    = [\sum_{m\in A, n}b_{mn}h_m\otimes t^{2n} +\sum_{u,v}c_{uv}x_u\otimes t^{2v+1}, \cr
   &   \,\,\,\,\,\,\,\,\,\, \frac{1}{2}j^{-1}h_i'\otimes t^{-2j}-
        \sum_{k\in B}d_kh_k'\otimes t^{-2j}+\sum_{p,q}e_{pq}x_p\otimes t^{2q+1}] \cr
   & = b_{ij}K+\sum_{m\in A, n,q}b_{mn}e_{mq}x_m\otimes t^{2n+2q+1}
          \,\,\,\,\, -\sum_{u,v}\beta_u(h_i')c_{uv}\frac{1}{2}j^{-1}x_u\otimes t^{2v-2j+1} \cr
      & \,\,\,\,\, +\sum_{m\in B, u,v}\beta_u(h_m')c_{uv}d_m x_u\otimes t^{2v-2j+1} \cr
      & = b_{ij}K+\sum_{k\in A, n,q}b_{k,n-q}e_{kq}x_k\otimes t^{2n+1}
       -\sum_{k,n}\beta_k(h_i')c_{k,n+j}\frac{1}{2}j^{-1}x_k\otimes t^{2n+1} \cr
   & \,\,\,\,\, +\sum_{m\in B,k,n}\beta_k(h_m')c_{k,n+j}d_m x_k\otimes t^{2n+1} \cr
   & = b_{ij}K+\sum_{k\in A,n}\{\sum_q b_{k,n-q}e_{kq}-\beta_k(h_i')c_{k,n+j}\frac{1}{2}j^{-1}\}x_k\otimes t^{2n+1} \cr
   &  \,\,\,\,\,+\sum_{k\in A,n}\sum_{m\in B}\beta_k(h_m')c_{k,n+j}d_mx_k\otimes t^{2n+1} \cr
   &  \, \,\,\,-\sum_{k\in B, n}\{\beta_k(h_i')c_{k,n+j}\frac{1}{2}j^{-1}- \sum_{m\in B}\beta_k(h_m')c_{k,n+j}d_m\}x_k\otimes t^{2n+1}.\cr} $$

For $k\in A$ and any $n$, the following single equation has solutions $e_{kq}$ with arbitrary values of $d_m$,
since its coefficients $b_{k, n-q}$ are not all zero
$$ \sum_q b_{k,n-q}e_{kq}-\beta_k(h_i')c_{k,n+j}\frac{1}{2}j^{-1}+\sum_{m\in B}\beta_k(h_m')c_{k,n+j}d_m=0.$$

For $k\in B$ and any $n$, the following system of equations on unknowns $d_m$
$$ \beta_k(h_i')c_{k,n+j}\frac{1}{2}j^{-1}=\sum_{m\in B}\beta_k(h_m')c_{k,n+j}d_m $$
or equivalently
$$ \beta_k(h_i')\frac{1}{2}j^{-1}=\sum_{m\in B}\beta_k(h_m')d_m, \,\,k\in B$$
has a unique solution, since the matrix $(\beta_k(h_m'))=(\beta_k, \beta_m)$ of its coefficients is non-degenerate.
The proof of the lemma is finished.
\hfill $ $ \qed
\enddemo

Now we state and prove the main result of this section, it implies that homogenous almost inner derivations of $\hat L$ are
determined completely.

\proclaim{Theorem 4.3} Let $\hat L$ be the twisted affinization of a minimal Q-graded subalgebra $L$ with
$\dim L=2\dim [L,L]$, then any odd almost inner derivation of $\hat L$ is an inner derivation, and any even almost inner derivation $D$ of $\hat L$
can be written uniquely as an infinite sum as follows
$$ D=\sum_{i,j}a_{ij}D_{ij}+ad y, \,\, a_{ij}\in \Bbb{F},\,\, y\in \hat L,$$
where the almost inner derivation $D_{ij}$ is defined as in Lemma 4.2, and $y$ is an even element of the $\Bbb{Z}_2$-graded Lie algebra $\hat L$.

Moreover, the infinite subset $\{D_{ij}+Inn(\hat L); \,\, 1\leq i\leq l, j\in \Bbb{Z}\}$ of the quotient space
$AID(\hat L)/Inn(\hat L)$ is linearly independent.
\endproclaim
\demo{Proof} According to Lemma 4.1, there is an element $y\in \hat L$ such that
the derivation $D-ad y$ satisfies the following relations:
$$ (D-ad y)(I\otimes Rt)=0, \,\, (D-ad y)(K)=0, $$
$$ (D-ad y)(H\otimes R)\subset \Bbb{F}K. $$

Set $a_{ij}K=(D-ad y)(h_i\otimes t^{2j}), 1\leq i\leq l, j\in \Bbb{Z}$, then $D-ad y=\sum_{i,j}a_{ij}D_{ij}$,
here $D_{ij}$ is the even derivation defined in Lemma 4.2.
If $D$ is odd and nonzero, then $D-ad y=0$, and $D$ is an inner derivation of $\hat L$; If $D$ is even, then the element $y$ is even,
and we obtain the infinite sum of the theorem as required.

Finally, we need to show that any even inner derivation $D$ of $\hat L$ is zero if it satisfies the relations \,(with $y=0$) of Lemma 4.1.
Suppose that the element
$Y=\sum_{u,v}d_{uv}h_u\otimes t^{2v}+\sum_{p,q}e_{pq}x_p\otimes t^{2q+1}\in \hat{L}$ satisfies $D(X)=[X,Y]$ for any $X\in \hat{L}$.
Let $D(h_i\otimes t^{2j})=a_{ij}K$ for some $a_{ij}\in \Bbb{F}$, we have the following identities
$$ \eqalign{  a_{ij}K  & = D(h_i\otimes t^{2j}+x_m\otimes t^{2n+1}) \cr
  & =[h_i\otimes t^{2j}+x_m\otimes t^{2n+1}, \sum_{u,v}d_{uv}h_u\otimes t^{2v}+\sum_{p,q}e_{pq}x_p\otimes t^{2q+1}] \cr
  & =\sum_u d_{u,-j}(h_i,h_u)2jK+\sum_qe_{iq}x_i\otimes t^{2j+2q+1}-\sum_{v}d_{mv}x_m\otimes t^{2n+2v+1}. \cr} $$

Set $m\neq i$, the above equation shows that $e_{iq}=0, \forall q$ and $d_{mv}=0, \forall v$.
Varying $i$ and $m$ if necessary, we deduce that all the elements $e_{pq}$ and $d_{uv}$ are zero,
and this implies that the derivation $D$ is zero.

The proof of the theorem is finished.
\hfill $ $ \qed
\enddemo

\vskip.3cm \Refs\nofrills{\bf REFERENCES}
\bigskip
\parindent=0.45in

\leftitem{[1]} D.Burde, K.Dekimpe, B.Verbeke, Almost inner derivations of Lie algebras,
 J. of Algebras and Its Applications 17(2018), no.11, 26 pages.

\leftitem{[2]} D.Burde, K.Dekimpe, B.Verbeke, Almost inner derivations of Lie algebras II,
ArXiv: 1905.08145(2019).

\leftitem{[3]} D.Burde, K.Dekimpe, B.Verbeke, Almost inner derivations of 2-step nilpotent Lie algebras of genus 2,
ArXiv: 2004.10567v1 [math.RA] 22 Apr 2020.

\leftitem{[4]} H.Dietrich, W.A.de Graaf, A computational approach to almost inner derivations,
ArXiv: 2403.05905(2024) [math.RA].

\leftitem{[5]} Yaxin Shen, Xiandong Wang, Almost inner derivations of affinizations of minimal Q-graded subalgebras, ArXiv: 2507.09484(2025) [math.RT].

\leftitem{[6]} I.Kenji, K.Yoshiyuki, Representation theory of the Virasoro Algebra,
Springer Monographs in Mathematics, Springer, 2011.

\leftitem{[7]} J. Milnor, Eigenvalues of the Laplace operator on certain manifolds,
Proc.Nat.Acad.Sci. U.S.A.51(1964), 542ff.

\leftitem{[8]} C.S.Gordon, E.N.Wilson, Isospectral deformations of compact solvmanifolds,
 J.Differen- tial Geom.19(1984), no.1, 214-256.

\leftitem{[9]} J.Humphreys, Introduction to Lie algebras and representation theory, GTM9, Springer, Fifth edition, 1987.

\leftitem{[10]} N.Jacobson, Lie algebras, Dover, New York, 1979.

\leftitem{[11]} R.Farnsteiner, Derivations and central extensions of finitely generated graded Lie algebras,
 J. of Algebra 118, 33-45(1988).

\leftitem{[12]} S.Azam, Derivations of tensor product of algebras, Comm. Algebra 36:(3)(2008), 905-927.

\leftitem{[13]} R.V.Kadison, Local derivations, J. Algebra 130, 494-509(1990).

\leftitem{[14]} Ayupov Sh.A., Kudaybergenov K.K., Local derivations on finite dimensional Lie algebras,
 Linear algebra and appl., Vol.493, 381-388(2016).

\leftitem{[15]} J.K.Adashev and T.K.Kurbanbaev, Almost inner derivations of some nilpotent Leibniz algebras,
 arXiv:2010.01904v2[math.RA] (2020).

\endRefs
\vfill
\enddocument
\end